\documentclass[a4paper,11pt]{article}
\usepackage{amsmath}
\usepackage{amsfonts}
\usepackage{amssymb}
\usepackage{enumerate}
\usepackage{graphics}
\RequirePackage{geometry}
\newtheorem{thm}{Theorem}
\newtheorem{lem}[thm]{Lemma}

\newtheorem{prop}[thm]{Proposition}
\newtheorem{defi}[thm]{Definition}

\numberwithin{equation}{section}

\newcommand{\prob}{\mathbb{P}}
\newcommand{\Trian}{\mathbb{T}}
\newcommand{\HP}{\mathbb{H}}
\newcommand{\RR}{\mathbb{R}}
\newcommand{\CC}{\mathbb{C}}
\newcommand{\imag}{\mathbf{i}}
\newcommand{\II}{\mathbf{1}}
\newcommand{\lra}{\leftrightarrow}
\newcommand{\cl}{\mathcal{C}}

\newcommand{\qed}{\hfill $\square$}

\begin{document}
\title{Factorization Formulas for $2D$ Critical Percolation, Revisited}
\author{R.P. Conijn\footnote{VU University Amsterdam, email: R.P.Conijn@vu.nl}}
\date{}
\maketitle

\begin{abstract}
We consider critical site percolation on the triangular lattice in the upper half-plane.
Let $u_1,u_2$ be two sites on the boundary and $w$ a site in the interior.
It was predicted by Simmons, Kleban and Ziff (2007) that the ratio
$\mathbb{P}(nu_1\leftrightarrow\,nu_2\leftrightarrow\,nw)^{2}\,/\,\mathbb{P}(nu_1\leftrightarrow\,nu_2)\cdot\mathbb{P}(nu_1\,\leftrightarrow\,nw)\cdot\mathbb{P}(nu_2\leftrightarrow\,nw)$
converges to $K_F$ as $n\to\infty$,
where $x\leftrightarrow\,y$ denotes that $x$ and $y$ are in the same cluster, and $K_F$ is a constant.
Beliaev and Izyurov (2012) proved an analog of this in the scaling limit.
We prove, using their result and a generalized coupling argument, the earlier mentioned prediction.
Furthermore we prove a factorization formula for
$\mathbb{P}(nu_2\leftrightarrow\,[nu_1,nu_1+s];\,nw\leftrightarrow\,[nu_1,nu_1+s])$, where $s>0$.
\end{abstract}

%
\medskip\noindent
\begin{center}
{\small {\it 2010 Mathematics Subject Classification.} 60K35 (82B43). }
\end{center}

\medskip\noindent

\section{Introduction and Main results.}
We consider critical site percolation on the triangular lattice. See \cite{G99} for a general introduction and \cite{S11,W09} for
more recent progress in two dimensional percolation.
A lot of attention has been given to crossing probabilities and critical exponents,
which are believed to be universal. In particular it is believed that in the continuum limit of
many two dimensional critical percolation models, crossing probabilities are conformally invariant.
However this has only been proved for site percolation on the triangular lattice by Smirnov \cite{S01}.
Another interesting question is whether it is possible to examine the  higher order correlation functions.
These are the functions $\mathbb{E}[X_{v_1}X_{v_2}\cdots X_{v_n}]$, where $v_i$ is a vertex and
$X_{v_i}=\II\{ 0 \lra v_i \}$ is the indicator function of the event that $v_i$ is in the open cluster
of the origin. A possible approach to compute these correlation functions might be via
factorization formulas.

To state our main results we consider the hexagonal lattice, where every center of a hexagon is a site
of the triangular lattice $\Trian$ in the closure of the upper half-plane $\HP := \{z \in \CC: \Im z > 0\}$.
In this lattice two neighbouring sites $x,y \in \Trian$ have $|x-y| = 1$.
By $\prob_{\eta}$ we denote the probability measure of critical percolation on $\eta\Trian$, for $\eta > 0$.
Let $\eta > 0$ and let the random set $Q \subset \overline{\HP}$ be the union of all hexagons for which the center is open.
The points $u,v \in \overline{\HP}$ are connected if $u,v$ are in the same connected component of $Q$. We denote this by $u \lra v$.
Let, for $u \in \eta\Trian$, $\cl(u)$ denote the open cluster containing $u$. Let, for $A \subset \overline{\HP}$,
\[
 \cl(A) := \bigcup_{u \in A\cap \eta\Trian} \cl(u).
\]
Further we will denote the hypergeometric function by $\,_2F_{1}(a,b;c;d)$ (see for example
\cite{AS65Book}). We denote by $\mathbb{S} := \{ z \in \CC \,:\, \Im(z) \in (0,1), \Re(z) > 0 \}$ the semi-infinite strip.

Our first main result is a factorization formula for the probability that three given vertices are in the same cluster,
where two of the vertices are on the boundary of the half-plane.
\begin{thm}\label{thm:mainTheorem:factPPP}
 Let $u_{1}, u_{2} \in \RR$ and $w \in \HP$ and $u_{1} \neq u_{2}$, then
 \begin{equation}
  \lim_{\eta \to 0} \frac{\prob_{\eta}(u_{1}\lra u_{2} \lra w)^2}{\prob_{\eta}(u_{1} \lra u_{2})\prob_{\eta}(u_{1} \lra w)\prob_{\eta}(u_{2} \lra w)} = K_{F},
 \end{equation}
 where
 \[
  K_{F} = \frac{2^{7}\pi^{5}}{3^{3/2}\Gamma(1/3)^{9}}.
 \]
\end{thm}

This factorization formula was heuristically derived, using Conformal Field Theory arguments, by Simmons, Kleban and Ziff in \cite{SKZ07}.
Using the convergence of percolation exploration interfaces to $SLE_6$ (See e.g. \cite{S00,S01}),
a mathematical rigorous proof of an analog of this formula in the continuum scaling limit was given by Beliaev and Izyurov in \cite{BI12}.
See Theorem \ref{thm:factFormBI} for their result.
That result is the starting point in the proof of Theorem \ref{thm:mainTheorem:factPPP}.
To obtain Theorem \ref{thm:mainTheorem:factPPP} from it we state and prove a quite general and robust form of a coupling
result for one-arm like events (see Proposition \ref{lem:Coupling} in Section \ref{subsec:Coupling}).

Our second main result involves the limiting behaviour of the probability \newline
$\prob(\{u_2, w\} \subset \cl([u_1,u_1+s]))$,
where $u_1,u_2$ are on the boundary of the half-plane and $w$ is in the half-plane.
We have the following theorem.
\begin{thm}\label{thm:mainTheorem:factPPI}
 Let $u_{1} \in \RR, w\in \HP, s > 0$ and $u_2 > u_1+s$, then
 \begin{equation}\label{eq:thm:PPI}
  \lim_{\eta \to 0} \frac{\prob_{\eta}(\{u_2, w\} \subset \cl([u_1,u_1+s]))}{\prob_{\eta}(w \in \cl([u_1,u_1+s]) )\,\, \prob_{\eta}(u_2 \in \cl([u_1,u_1+s]))} = \psi(u_1,s,u_2,w),
 \end{equation}
 where $\psi$ is the function
 \begin{equation*}
 \psi(u_1,s,u_2,w) = e^{\pi x/3}\cdot \frac{\,_2F_1\left(-\frac{1}{2},-\frac{1}{3};\frac{7}{6};e^{-2\pi x}\right)}{\,_2F_1\left(-\frac{1}{2},-\frac{1}{3};\frac{7}{6};1\right)},
 \end{equation*}
with $x = \Re (\Psi_{u_1,s,u_2}(w))$ where
$\Psi_{u_1,s,u_2}$ is the conformal map that transforms \newline
$\{\HP, u_1,u_1+s,u_2\}$ to $\{ \mathbb{S},\imag,0,\infty\}$.
\end{thm}
Simmons, Ziff and Kleban studied in \cite{SZK09} the probability in the numerator in \eqref{eq:thm:PPI}.
They used Conformal Field Theory arguments to find several predictions for formulas of the probabilities
in \eqref{eq:thm:PPI}.
Theorem \ref{thm:mainTheorem:factPPI} is a discrete analog of one of their predictions
(Equation (29) in Section III B of \cite{SZK09}).

Our interest in these factorization formulas came from the paper \cite{BI12} by Beliaev and Izyurov.
They rigorously proved an analog of the formula \eqref{eq:thm:PPI} above in the scaling limit,
but with the probability $\prob(w \in \cl([u_1,u_1 + s]))$
replaced by $s_{3}^{5/48}$, see Theorem \ref{thm:condProbBI}.
However their theorem involves probabilities where the cluster does not necessarily touch $w$, but comes
within a certain distance from it.
More precisely, their formula is about the limits where first the mesh size, and secondly the above mentioned distance
tends to zero.

\medskip
\noindent\textbf{Remark: }We believe that our coupling argument, Proposition \ref{lem:Coupling},
is more generally applicable. For example Simmons, Ziff and Kleban also predicted in \cite{SZK09}
a factorization formula for the probability
$\prob_{\eta}(u_2 \lra w \lra [u_1,u_1+s])$.
We hope that as soon as an analog of this result in the scaling limit
has been proved, our Proposition \ref{lem:Coupling} can be used to prove this factorization formula
in a discrete setting.
More recently Delfino and Viti heuristically derived in \cite{DV11} (see also \cite{ZSK11}) a
factorization formula for the probability
$\prob(x\lra y\lra w)$, where all three points are in the interior of the half-plane.
We also believe that Proposition \ref{lem:Coupling} might be an ingredient for a rigorous proof of a discrete analog
of this factorization formula, again after the scaling limit analog has been proved.

\bigskip
The rest of the paper is organized as follows.
In Section \ref{sec:preliminaries} we introduce some notation and sum up some preliminary results,
which are crucial for our proofs. In Section \ref{subsec:Coupling} we state and proof a quite general and
abstract ratio limit result, Proposition \ref{lem:Coupling}, which is based on a coupling argument.
This proposition forms a key ingredient for the proofs of both main theorems.
In the last Sections \ref{subsec:pfOfThmfactPPP} and \ref{subsec:pfOfThmfactPPI} we give the proofs of
our main results.


\section{Notation and Preliminaries.}\label{sec:preliminaries}
We begin with some notation. Let $\Omega^{\eta} := \{0,1\}^{\eta\Trian}$.
Elements of $\Omega^{\eta}$ will typically be denoted by $\omega,\nu$ and called \textit{configurations}.
We call a vertex $v \in \eta\Trian$ \textit{open} if $\omega_{v} = 1$, otherwise we say that $v$ is \textit{closed}.
For two configurations $\omega,\nu \in \Omega^{\eta}$ we write $\omega \le \nu$ if and only if
$\omega_{v} \le \nu_{v}$ for all $v \in \eta\Trian$.
Let $P \subset \HP$, we write $\omega_{P} \in \{ 0,1 \}^{\eta\Trian \cap P}$ for the restriction of
$\omega$ to the vertices which are contained in $P$.
For two disjoint sets $P,Q \subset \HP$, and configurations $\omega_{P}, \omega_{Q}$
we define $\omega_{P} \times \omega_{Q}$ to be the configuration
$\tilde{\omega}_{P\cup Q} \in \{ 0,1 \}^{\eta\Trian \cap (P\cup Q)}$
such that $\tilde{\omega}_{P} = \omega_{P}$ and $\tilde{\omega}_{Q} = \omega_{Q}$.
Let $V \subset \Omega^{\eta}$ be an event and $A \subset \HP$. We define the event
\begin{equation}\label{eq:def:V_A}
 V_{A} := \{ \omega\,|\, \exists\, \tilde{\omega}_{\HP\setminus A}: \omega_{A}\times \tilde{\omega}_{\HP\setminus A} \in V\}.
\end{equation}
Further, with some abuse of notation, for $A \subset \HP, \omega_{A} \in \{0,1\}^{A\cap \eta\Trian}$ and
$V \subset \Omega^{\eta}$ we write $\prob_{\eta}(V \,|\, \omega_{A})$
for the conditional probability of $V$ given that the configuration on $A$ equals $\omega_{A}$.
Similarly we write $\{\omega_{A}\}$ for the event that the configuration on $A$ equals $\omega_{A}$.

For $z = z_1+z_2\imag \in \HP$ and $a > 0$, we write $B_{a}(z)$
for the intersection of the half-plane with the $2a \times 2a$-box centered at $z$.
We denote annuli by $A(z; a,b) := B_{b}(z) \setminus B_{a}(z)$.
A \textit{circuit} in an annulus $A(z;a,b)$ is a sequence of neighbouring vertices in $\eta\Trian$, such that every vertex appears at most once
in the sequence, the last vertex is a neighbour of the first and it surrounds $B_{a}(z)$. We will often encounter
annuli which intersect the boundary of $\HP$, in that case we will also consider \textit{semi-circuits}.
A semi-circuit in an annulus $A(z;a,b)$ is a sequence of neighbouring vertices such that every vertex appears at
most once in the sequence, the first and the last vertex are both on the boundary $\partial \HP$ and the
semi-circuit 'surrounds' $B_{a}(z)$. In other words a semi-circuit is a path in $\HP$ from the boundary of
$\HP$ to the boundary of $\HP$ which disconnects $B_{a}(z)$ from infinity.
A (semi-)circuit is called open if all its vertices are open.
For a (semi-)circuit $\gamma$ we denote by $int(\gamma)$ the bounded connected component of $\HP \setminus \bar{\gamma}$
containing $B_{a}(z)$, where $\bar{\gamma}$ is the curve in the plane described by $\gamma$. Further $ext(\gamma)$ is the
unbounded connected component of $\HP \setminus \bar{\gamma}$.


Let $\mathbb{U} := \{ z \in \CC: |z| < 1 \}$ be the open ball of radius one. For $w \in \HP$ and
a closed connected set $A \subset \HP$ we denote by
$\rho(w,A)$ the conformal radius of the component of $w$ in $\HP\setminus A$ seen from $w$.
It is defined as follows. If $w \not\in A$, let $V$ be the connected component of $w$ in $\HP \setminus A$. Let
$\phi: V \to \mathbb{U}$ be the unique conformal map with $\phi(w) = 0$ and $\phi'(w) > 0$.
Then we set $\rho(w,A) := 1/\phi'(w)$.
Otherwise, if $w \in A$ we set $\rho(w,A) := 0$.
We can compare the conformal radius with the euclidean distance from the point to the set, namely it follows from Koebe's 1/4-Theorem
and Schwarz' Lemma that
\begin{equation}\label{eq:pf:mainThm:equivConfRad}
 \frac{1}{4}\rho(w,A) \le \min_{x\in A} |w-x| \le \rho(w,A).
\end{equation}
(See e.g. \cite{A73Book})

\begin{figure}
 \centering
 \scalebox{1.0}{\includegraphics{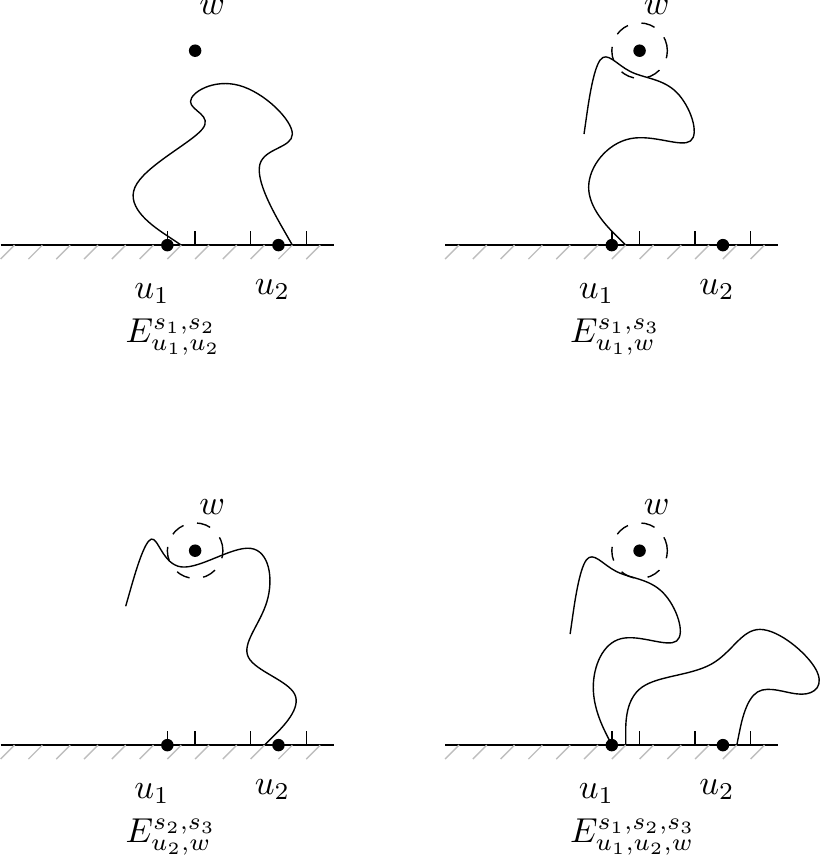}}
 \caption{The events $E_{u_1,u_2}^{s_1,s_2}, E_{u_1,w}^{s_1,s_3}, E_{u_2,w}^{s_2,s_3}$ and $E_{u_1,u_2,w}^{s_1,s_2,s_3}$. Note that the clusters in
 $E_{u_1,u_2,w}^{s_1,s_2,s_3}$ might be disjoint.}
 \label{fig:eventsE}
\end{figure}
We introduce the following events, which all represent the existence of clusters which come close to certain vertices. See Figure \ref{fig:eventsE}.
For $u_1,u_2 \in \RR$, $w \in \HP$ and $s_1,s_2,s_3 > 0$,
\begin{eqnarray}
 E_{u_1,u_2}^{s_1,s_2} & := & \{ \cl([u_1,u_1+s_1]) \cap [u_2 - s_2, u_2+s_2] \neq \emptyset\};  \label{eq:defi:E1}\\
 E_{u_1,w}^{s_1,s_3} & := & \{\rho(w,\cl([u_1,u_1+s_1])) < s_3\}; \nonumber\\
 E_{u_2,w}^{s_2,s_3} & := & \{\rho(w,\cl([u_2 - s_2, u_2+s_2])) < s_3\}; \nonumber\\
 E_{u_1,u_2,w}^{s_1,s_2,s_3} & := & E_{u_1,u_2}^{s_1,s_2} \cap E_{u_1,w}^{s_1,s_3}.\nonumber
\end{eqnarray}
Although all these events depend on $\eta$, we omit this from the notation.
They represent the discrete versions of the events used by Beliaev and Izyurov in \cite{BI12}.
Note the difference between the events $E_{u_1,w}^{s_1,s_3}$ and $E_{u_2,w}^{s_2,s_3}$. This is to stay as close
as possible to the events defined in that paper.
As mentioned before Beliaev and Izyurov considered the limits, as $\eta \to 0$, of the probabilities of the
events above. That is
\begin{eqnarray*}
 f_{u_1,u_2}^{s_1,s_2} & := & \lim_{\eta \to 0} \prob_{\eta}(E_{u_1,u_2}^{s_1,s_2}); \\
 f_{u_1,w}^{s_1,s_3} & := & \lim_{\eta \to 0} \prob_{\eta}(E_{u_1,w}^{s_1,s_3}); \\
 f_{u_2,w}^{s_2,s_3} & := & \lim_{\eta \to 0} \prob_{\eta}(E_{u_2,w}^{s_2,s_3}); \\
 f_{u_1,u_2,w}^{s_1,s_2,s_3} & := & \lim_{\eta \to 0} \prob_{\eta}(E_{u_1,u_2,w}^{s_1,s_2,s_3}).
\end{eqnarray*}
The existence of these limits follow from the results in \cite{LSW02,S01}.
Namely the existence of the first one (which is actually given by Cardy's formula) was proved by Smirnov in \cite{S01}.
The second and third are described in the article on the one-arm exponent for critical $2D$ percolation \cite{LSW02},
using the so called exploration path, started at, respectively $u_1+s_1$ and $u_2 + s_2$.
The fourth one can also be described in terms of exploration path.
It is the intersection of the events: (1) the exploration path starting at $u_1+s_1$ swallows $u_2-s_2$ before it swallows
$u_1$ or $u_2 + s_2$ and (2) the exploration path, or union of nested exploration paths,
comes $s_3$ close to $w$ in conformal radius. See \cite{LSW02} for the definition of the exploration path and more details.

As Beliaev and Izyurov already mentioned in \cite[Remark 4]{BI12}, the factorization formula they
proved, Proposition 4.1 in their paper, implies the following Theorem.
\begin{thm}[Remark 4 in \cite{BI12}]\label{thm:factFormBI}
 Let $u_{1}, u_{2}, w$ and $K_F$ be as in Theorem \ref{thm:mainTheorem:factPPP}.
 For every $\varepsilon, s_0 > 0$ there exist $s_1,s_2,s_3 < s_0$ such that
 \begin{equation}
  \left|\frac{(f_{u_{1}, u_{2}, w}^{s_1,s_2,s_3})^2}{f_{u_{1}, u_{2}}^{s_1,s_2}\cdot f_{u_{1}, w}^{s_1,s_3}\cdot f_{u_{2}, w}^{s_2,s_3}} - K_{F}\right| < \varepsilon.
 \end{equation}
\end{thm}
The following Theorem is the main result in \cite{BI12}, and will be used in the proof of Theorem \ref{thm:mainTheorem:factPPI}.
\begin{thm}[Theorem 1.1 in \cite{BI12}]\label{thm:condProbBI}
 Let $u_1,u_2,w,s$ be as in Theorem \ref{thm:mainTheorem:factPPI}.
 One has
 \begin{equation}
  \lim_{s_3 \to 0}\lim_{s_2 \to 0}s_{3}^{-5/48}\cdot \frac{f_{u_1,u_2,w}^{s,s_2,s_3}}{f_{u_1,u_2}^{s,s_2}} = K_1|\Psi'_{u_1,s,u_2}(w)|^{5/48}G\left(\Re (\Psi_{u_1,s,u_2}(w)), \Im (\Psi_{u_1,s,u_2}(w))\right),
 \end{equation}
 where $\Psi_{u_1,s,u_2}$ is the conformal map that transforms $\{\HP, u_1,u_1+s,u_2\}$ to $\{ \mathbb{S},\imag,0,\infty\}$
 and
 \begin{eqnarray}
  K_1 & = & \frac{18 \pi^{5/48}}{5\pi\cdot 2^{5/48}} H(0)^{-1}\nonumber\\
  G(x,y) & = & e^{\pi x/3} H(x)\sinh(\pi x)^{-1/3} \left(\frac{\sinh(\pi x)^{2}\sin(\pi y)^{2}}{\sinh(\pi x)^{2}+\sin(\pi y)^{2}} \right)^{11/96},\label{eq:thmCondProbBI:DefG}
 \end{eqnarray}
 with
 \begin{equation}
  H(x) = \,_2F_1\left(-\frac{1}{2},-\frac{1}{3};\frac{7}{6};e^{-2\pi x}\right).
 \end{equation}
\end{thm}
The lemma below, proved by Beliaev and Izyurov, is an improvement of a result by Lawler, Schramm and Werner in \cite{LSW02}.
\begin{lem}[Lemma 2.2 in \cite{BI12}]\label{lem:poiInt}
 Let $u_{1}, w$ be as in Theorem \ref{thm:mainTheorem:factPPP} and let $s > 0$.
 One has
 \begin{equation}
  \lim_{s_3 \to 0} s_3^{-5/48}\cdot f_{u_1,w}^{s,s_3} = K_2 |\phi'(w)|^{5/48}(\sin(\pi\omega/2))^{1/3},
 \end{equation}
 where $\omega$ is the harmonic measure of $(u_1,u_1 + s)$ seen from $w$; $\phi$ is a conformal map from $\HP$ to the unit disc such that $\phi(w) = 0$, and
 \begin{equation}
  K_2 = \frac{18}{5\pi}.
 \end{equation}
\end{lem}

We end this section with a lemma which is a simple generalization of the FKG inequality.
\begin{lem}\label{lem:genFKG}
 Let $A \subset \HP$ and let $B,E$ be increasing events. Let $\nu_{A} \in \{ 0,1 \}^{\eta\Trian \cap A}$.
 If $B$ is completely determined by the vertices in $\HP \setminus A$, that is $B=B_{\HP\setminus A}$,
 then
 \[
  \prob_{\eta}(B\cap E\cap \{\nu_{A}\} ) \ge \prob_{\eta}(B)\prob_{\eta}(E \cap \{\nu_{A}\}).
 \]
\end{lem}

\noindent\textit{Proof of Lemma \ref{lem:genFKG}:}
The proof of this lemma is straightforward and we omit it.
\qed

\section{Proofs of the main results.}
%

\subsection{Coupling of one-arm like events.}\label{subsec:Coupling}
The proof of our first main result, Theorem \ref{thm:mainTheorem:factPPP}, has two ingredients. The first is Theorem \ref{thm:factFormBI}.
The second ingredient for our proof is a coupling argument for one-arm like events
which appeared in somewhat different forms in \cite{K86} and more recently in \cite{GPS13}. However our
coupling result is developed in a more general framework of one-arm like events;
see Definitions \ref{defi:OneArmLike}-\ref{defi:comparable} below.

Our second main result, Theorem \ref{thm:mainTheorem:factPPI},
also has this coupling argument as one of the main ingredients. The other main ingredients for the proof of Theorem
\ref{thm:mainTheorem:factPPI} are Theorem \ref{thm:condProbBI} and Lemma \ref{lem:poiInt}.

The proof of our coupling argument is along the lines of the sketch in \cite{GPS13}. In that paper,
among other very interesting results,
a ratio limit theorem was proved. They proved that, for every $a>0$
\begin{equation*}
 \lim_{\eta \to 0}\, \frac{\prob_{\eta}(0 \lra \CC\setminus [-a,a]^2)}{\prob_{\eta}(0 \lra \CC\setminus [-1,1]^2)} = a^{-5/48},
\end{equation*}
see section 5.1 in that paper.
Here we show that their arguments can be modified, which makes them more generally applicable.
In the arguments of \cite{GPS13}, when a cluster comes $s$ close to a point $z$ it means that
the cluster touches the boundary of $B_{s}(z)$. Hence the configuration in $B_{s}(z)$ is 
independent of the event that the cluster comes close. However, in our situation,
when a cluster comes close to a vertex $z$ it means in some occasions that the conformal radius is small
and in other occassions it means that the cluster touches the interval $[z-s,z+s]$, as we saw in Section \ref{sec:preliminaries}.
Hence in our situation the configuration in $B_{s}(z)$ is not independent from the event that the cluster
comes $s$ close to $z$.
This difference in measuring the distance of a cluster to a point makes the arguments
more complicated.
Our way to solve these complications is to grasp the essence which makes things work.
This led us to the following formal definition of a class of events
which intuitively describe the occurrence of a cluster coming within a distance $s$ from $z$.
\begin{defi}\label{defi:OneArmLike}
 Let $s, C > 0$.
 Let $z \in \HP$ and $V \subset \Omega^{\eta}$ be an increasing event.
 We say that $V$ is an \emph{$(s,C)$-one-arm like event around $z$}
 if, for every (semi-)circuit $\gamma$ in $A(z;s,C)$,
 \begin{equation}
  V \left\{\begin{array}{cl}
                    \subset & \{ B_{s}(z) \lra \HP \setminus B_{C}(z) \}\\
                    \supset &  \{ \gamma \textrm{ open}\} \cap V_{ext(\gamma)} \cap V_{int(\gamma)}
                   \end{array}\right.
 \end{equation}
 and
 \[
  \{I(z,s) \lra \gamma\} \subset V_{int(\gamma)},
 \]
 where $I(z,s)$ is the horizontal line segment $[z,z+s/8] \subset \overline{\HP}$ and
 $V_{int(\gamma)}, V_{ext(\gamma)}$ as in \eqref{eq:def:V_A}.
\end{defi}
For example, for every $x,s,C \in \RR$ and $a \in [1/8,1]$, the events $\{B_{as}(x\imag) \lra (x\imag+2C(1+\imag)))\}$
and $\{I(x,s) \lra \HP \setminus B_{2C}(x)\}$ are $(s,C)$-one-arm like events around $x\imag$, respectively $x$.
In the proof of Theorem \ref{thm:mainTheorem:factPPP} we will see that also certain events
concerning a small conformal radius from $z$ to a certain cluster are $(s,C)$-one-arm like events.

Observe that the definition above implies that for every (semi-)circuit $\gamma$ in $A(z;s,C)$,
\begin{equation}\label{eq:obs:OneArmIndep}
 V \cap \{ \gamma \textrm{ open} \} = V_{ext(\gamma)} \cap V_{int(\gamma)} \cap \{ \gamma \textrm{ open} \},
\end{equation}
where $V$ is an $(s,C)$-one-arm like event around $z$.

If $V$ is an $(s,C)$-one-arm like event around $z$,
there is a certain open cluster which comes within a distance $s$ from $z$.
For any such event $V$ we will also consider a related event where this cluster hits $z$.
Intuitively a good candidate for such an event would be $V \cap \{ z \lra \HP \setminus B_{C}(z)\}$,
but this is not appropriate: under this event the cluster $\cl(z)$ and the earlier mentioned cluster,
could be disjoint. In other words, this event is too large.
It turns out that the following definition is suitable for our purposes.

\begin{defi}\label{defi:pointVersion}
 Let $V$ be an $(s,C)$-one-arm like event around $z$.
 Let $V^{\bullet}$ be an increasing event.
 We call $V^{\bullet}$ a \emph{point version of $V$}
 if, for every (semi-)circuit $\gamma$ in $A(z;s,C)$,
 \begin{equation}
  V^{\bullet} \left\{\begin{array}{cl}
                    \subset & V \cap \{ z \lra \HP \setminus B_{C}(z) \}\\
                    \supset &  \{ \gamma \textrm{ open}\} \cap V_{ext(\gamma)} \cap \{z \lra \gamma\}.
                   \end{array}\right.
 \end{equation}
\end{defi}
For example, for every $x,s,C \in \RR$ and $a \in [1/8,1]$, the event $\{x\imag \lra (x\imag+2C(1+\imag)))\}$ is
a point version of $\{B_{as}(x\imag) \lra (x\imag+2C(1+\imag)))\}$
and $\{x \lra \HP \setminus B_{2C}(x)\}$ is a point version of $\{I(x,s) \lra \HP \setminus B_{2C}(x)\}$.
To state the coupling proposition we need one more definition.
\begin{defi}\label{defi:comparable}
 Let $z \in \HP$ and $s,C>0$.
 Let $V$ and $W$ be $(s,C)$-one-arm like events around $z$.
 We say that $V, W$ are \emph{$(s,C)$-comparable around $z$} if
 the events $V_{B_{C}(z)}$ and $W_{B_{C}(z)}$ are equal.
\end{defi}
It follows easily from this definition, that equality also holds for any subset of $B_{C}(z)$.
In other words, let $V,W$ be $(s,C)$-comparable around $z$, then $V_{A} = W_{A}$ for every $A \subset B_{C}(z)$.

Our coupling argument is contained in the following proposition.
\begin{prop}\label{lem:Coupling}
 Let $C>0$ and $z \in \HP$. There exist increasing functions $\varepsilon(s), m(s): \RR_{+} \to (0,1)$,
 with $\varepsilon(s) \to 0$ and $m(s) \to 0$ as $s \to 0$
 such that the following holds. For all $s>0$, for all $\eta < m(s)$ and for every pair
 $V, W \subset \Omega^{\eta}$ of $(s,C)$-comparable events around $z$ and
 point versions $V^{\bullet}$ of $V$ and $W^{\bullet}$ of $W$ we have
 \begin{equation}\label{eq:lem:Coupling}
  \left| \frac{\prob_{\eta}(V^{\bullet}\,|\, V)}{\prob_{\eta}(W^{\bullet}\,|\, W)} - 1 \right| < \varepsilon(s).
 \end{equation}
\end{prop}

Before we give a proof of this proposition, we introduce some notation and state a lemma
which is crucial in the proof of Proposition \ref{lem:Coupling}.

Let $C,s > 0$ and $z \in \HP$. Let $l(i) := 4^{-i}C$. Let $N(s,C) = \lfloor\log_{4}(C / s)\rfloor - 2$ and
let $P_i := \HP \setminus B_{l(i)}(z)$.
We define for every $i \in \{0, 1,2, \cdots, N(s,C)\}$
the annuli $AI_{i} := A(z;\frac{1}{4}l(i),\frac{1}{2}l(i))$,
$AO_{i} := A(z;\frac{1}{2}l(i),l(i))$ and $A_{i} := AI_{i} \cup AO_{i}$.
We denote by $\Gamma I_{i}$ the outermost open (semi-)circuit
in $AI_{i}$ and by $\Gamma O_{i}$ the innermost open (semi-)circuit in $AO_{i}$, if they exist.
Otherwise, if there is no (semi-)circuit in $AI_{i}$ (resp. $AO_{i}$)
we set $\Gamma I_{i} = \emptyset$ (resp. $\Gamma O_{i} = \emptyset$).
Let $\gamma_{I}$ be a fixed (semi-)circuit in $AI_{i}$ and $\gamma_{O}$ be a
fixed (semi-)circuit  in $AO_{i}$.
The following observation is quite standard.
Conditioned on $\{\Gamma I_{i} = \gamma_{I}; \Gamma O_{i} = \gamma_{O}\}$, the configuration in
$int(\gamma_{I}) \cup ext(\gamma_{O})$ is a fresh independent copy of a percolation configuration.
\begin{lem}\label{lem:CouplKeyIngr}
 There exists a universal constant $C_1 \in (0,1)$ such that the following holds.
 Let $z \in \HP, s,C > 0, i \le N(s,C)$ and let $\gamma_{I}$ be a deterministic (semi-)circuit.
 Let $V$ be an $(s,C)$-one-arm like event around $z$.
 Then, for every $\nu \in V_{P_i}$ we have
 \begin{equation}\label{eq:lem:couplCompProb}
 \prob_{\eta}(\Gamma I_{i} = \gamma_{I} \,|\, V\cap \{\nu_{P_{i}}\})
  \ge C_1\, \prob_{\eta}(\{\Gamma I_{i} = \gamma_{I}\} \cap \{\Gamma O_{i}\textrm{ exists}\} \cap \{\gamma_{I} \lra \Gamma O_{i}\}).
\end{equation}
\end{lem}

\noindent\textit{Proof of Lemma \ref{lem:CouplKeyIngr}:} 
It is sufficient to prove that, for every (semi-)circuit $\gamma_{O}$,
\begin{eqnarray}
  \lefteqn{\prob_{\eta}(\{\Gamma I_{i} = \gamma_{I}\} \cap \{\Gamma O_{i} = \gamma_{O}\}\cap \{\gamma_{I} \lra \gamma_{O}\} \,|\, V\cap \{\nu_{P_{i}}\})}\label{eq:pf:lem:couplCompProb}\\
& \ge & C_1\, \prob_{\eta}(\{\Gamma I_{i} = \gamma_{I}\}\cap \{\Gamma O_{i} = \gamma_{O}\}\cap \{\gamma_{I} \lra \gamma_{O}\}).\nonumber
 \end{eqnarray}
Namely \eqref{eq:lem:couplCompProb} immediately follows from \eqref{eq:pf:lem:couplCompProb} after summing
over the possible (semi-)circuits $\gamma_{O}$.

Let $\gamma_{O}$ be an arbitrary (semi-)circuit and
\[
 D=\{\Gamma I_{i} = \gamma_{I}\}\cap \{ \Gamma O_{i} = \gamma_{O}\}\cap \{\gamma_{I} \lra \gamma_{O}\}.
\]
Then
the left hand side of \eqref{eq:pf:lem:couplCompProb} is equal to
\begin{equation}\label{eq:pf:CouplIngr:condProb}
 \frac{\prob_{\eta}(D\cap V\cap \{\nu_{P_{i}}\})}{\prob_{\eta}(V\cap \{\nu_{P_{i}}\})}.
\end{equation}
It follows from \eqref{eq:obs:OneArmIndep} and Definition \ref{defi:OneArmLike} that
\begin{eqnarray*}
 \prob_{\eta}(D \cap V\cap \{\nu_{P_{i}}\}) & = & \prob_{\eta}(D\cap V_{ext(\gamma_{O})}\cap V_{int(\gamma_{O})}\cap \{\nu_{P_{i}}\})\\
 & \ge & \prob_{\eta}(D\cap V_{ext(\gamma_{O})}\cap \{I(z,s) \lra \gamma_{I}\}\cap \{\nu_{P_{i}}\}).
\end{eqnarray*}
The last probability is, by the observation about inner- and outermost (semi-)circuits, equal to
\begin{equation}\label{eq:pf:CouplIngr:teller}
 \prob_{\eta}(D)\prob_{\eta}(I(z,s) \lra \gamma_{I})\prob_{\eta}(V_{ext(\gamma_{O})}\cap \{\nu_{P_{i}}\}).
\end{equation}
On the other hand the denominator in \eqref{eq:pf:CouplIngr:condProb} is, again by Definition \ref{defi:OneArmLike}, less than or equal to
\begin{eqnarray}\label{eq:pf:CouplIngr:noemer}
 \prob_{\eta}(V_{ext(\gamma_{O})}\cap \{\nu_{P_i}\}\cap \{B_{s}(z) \lra \gamma_{I}\}) & = & \prob_{\eta}(V_{ext(\gamma_{O})}\cap \{\nu_{P_i}\})\prob_{\eta}(B_{s}(z) \lra \gamma_{I}) \\
  & \le & \prob_{\eta}(V_{ext(\gamma_{O})}\cap \{\nu_{P_i}\})\cdot \frac{1}{C_{1}}\prob_{\eta}(I(z,s) \lra \gamma_{I}),\nonumber
\end{eqnarray}
where the constant $C_1$ comes from standard RSW and FKG arguments.
A combination of \eqref{eq:pf:CouplIngr:condProb}, \eqref{eq:pf:CouplIngr:teller} and \eqref{eq:pf:CouplIngr:noemer}
gives \eqref{eq:pf:lem:couplCompProb}. This finishes the proof of Lemma \ref{lem:CouplKeyIngr}.
\qed

\medskip\noindent\textit{Proof of Proposition \ref{lem:Coupling}:}
We will describe a coupling of the conditional distributions given $V$ and given $W$, denoted by $\tilde{\prob}$.
More precisely we construct $\tilde{\prob}$ such that, for $\nu,\omega \in \Omega^{\eta}$,
\begin{equation}\label{eq:pf:Coupl:CouplMeas}
 \tilde{\prob}(\nu \times \Omega^{\eta}) = \prob_{\eta}(\nu \,|\, V), \qquad \tilde{\prob}( \Omega^{\eta} \times \omega) = \prob_{\eta}(\omega \,|\, W).
\end{equation}
Furthermore $\tilde{\prob}$ will be such that the probability that the two distributions are
successfully coupled (in a sense defined precisely below) goes to 1 as $s$ tends to zero, uniformly in $\eta$.
We will finish the proof by showing how this coupling can be used to prove the proposition.

Let us first describe the coupling procedure.
First we draw, independently of each other, $\nu_{P_{0}}$
and $\omega_{P_{0}}$ according to, respectively $\prob_{\eta}(\cdot \,|\, V)$
and $\prob_{\eta}(\cdot \,|\, W)$.
Next we draw, step by step, the random elements $\nu_{A_{i}}$, $\omega_{A_{i}}$,
starting from $i=0$.

Every step goes as follows. The outermost (semi-)circuits $\Gamma I_{i}(\nu)$,
$\Gamma I_{i}(\omega)$ are drawn from the optimal coupling
of $\prob_{\eta}(\Gamma I_{i}(\nu) = \cdot \,|\, V; \nu_{P_{i}})$
and $\prob_{\eta}(\Gamma I_{i}(\omega) = \cdot \,|\, W; \omega_{P_{i}})$.
That is, the coupling is such that
$\tilde{\prob}(\Gamma I_{i}(\nu) = \Gamma I_{i}(\omega) \neq \emptyset \,|\, \nu_{P_{i}};\, \omega_{P_{i}})$
is as large as possible.

We say that this step of the coupling is successful if $\Gamma I_{i}(\nu) \neq \emptyset$ and
$\Gamma I_{i}(\nu) = \Gamma I_{i}(\omega) =: \gamma$.
In that case we can finish the coupling procedure as follows.
First we draw $\nu_{ext(\Gamma I_{i}(\nu))\cap A_{i}}$ and $\omega_{ext(\Gamma I_{i}(\omega))\cap A_{i}}$
from the appropriate conditional probability measures, independently of each other.
So $\nu_{ext(\Gamma I_{i}(\nu))\cap A_{i}}$ is drawn from the probability measure
$\prob_{\eta}(\cdot \,|\, \Gamma I_{i}(\nu) = \gamma; V; \nu_{P_{i}})$.
Since $V$ is an $(s,C)$-one-arm like event we have for every
$\nu_{int(\gamma)} \in \{0,1\}^{\eta\Trian \cap int(\gamma)}$
\begin{eqnarray*}
 \prob_{\eta}(\nu_{int(\gamma)} \,|\, \Gamma I_{i}(\nu) = \gamma; V; \nu_{ext(\gamma)}) & = & \prob_{\eta}(\nu_{int(\gamma)} \,|\, V_{int(\gamma)};\, V_{ext(\gamma)};\, \Gamma I_{i}(\nu) = \gamma;\, \nu_{ext(\gamma)})\\
 & = & \prob_{\eta}(\nu_{int(\gamma)} \,|\, V_{int(\gamma)}),
\end{eqnarray*}
where we used \eqref{eq:obs:OneArmIndep} in the first equality and
independence of $\nu_{int(\gamma)}$ and $V_{int(\gamma)}$ from the rest in the second.
The same holds for $W$. Now we use that $V$ and $W$ are $(s,C)$-comparable around $z$.
As we saw immediately after Definition \ref{defi:comparable} this implies that $V_{int(\gamma)} = W_{int(\gamma)}$, hence
the two conditional distributions of the interior of $\gamma$ are equal.
Thus we can draw $\nu_{int(\gamma)}$
according to $\prob_{\eta}(\cdot \,|\, V_{int(\gamma)})$
and take $\omega_{int(\gamma)} := \nu_{int(\gamma)}$.

If this step of the coupling was not successful,
let $\gamma_{\nu}$ and $\gamma_{\omega}$ be the outcome of $\Gamma I_{i}(\nu)$ and
$\Gamma I_{i}(\omega)$ respectively,
we draw the random elements $\nu_{A_{i}}$, $\omega_{A_{i}}$
according to\newline $\prob_{\eta}(\cdot \,|\, \Gamma I_{i}(\nu) = \gamma_{\nu}; V; \nu_{P_{i}})$
and $\prob_{\eta}(\cdot \,|\, \Gamma I_{i}(\omega) = \gamma_{\omega}; W; \omega_{P_{i}})$
independently of each other and continue to the next step with $i+1$.


If all steps, $i=0,\cdots, N(s,C)$, of the coupling were not successful,
we draw $\nu_{RM}$ and $\omega_{RM}$
according to the appropriate conditional probabilities, independently of each other, where
\begin{equation}\label{eq:pf:Coupl:defRM}
RM := B_{l(N(s,C)+1)}(z) \supset B_{2s}(z).
\end{equation}

That this procedure defines a coupling for the measures in \eqref{eq:pf:Coupl:CouplMeas}
follows from standard arguments.

Let $S$ denote the event that the coupling is successful
(i.e. that some step in the above described procedure is succesful).
The crucial property of this coupling is that
\begin{equation}\label{eq:pf:Coupl:whenCoupled}
 (\Omega^{\eta} \times W^{\bullet})\cap S = (V^{\bullet} \times \Omega^{\eta})\cap S,
\end{equation}
which follows easily from Definition \ref{defi:pointVersion}.
To see that $\tilde{\prob}(S) \to 1$ as $s \to 0$, note that it follows easily from
Lemma \ref{lem:CouplKeyIngr} together with RSW, FKG arguments
that there exists a constant $C_{2} > 0$ such that for every $i$
\begin{eqnarray*}
 \sum_{\gamma_{I}} \min_{{E \in \{ V, W \}\atop \omega_{P_{i}} \in \{0,1\}^{P_{i}}}}\left(\prob_{\eta}(\Gamma I_{i} = \gamma_{I} \,|\, E;\, \omega_{P_{i}})\right) \ge C_2.
\end{eqnarray*}
Hence, for every step in the procedure described above,
the probability that the coupling is successful is at least $C_{2}$.
Thus
\begin{equation}\label{eq:pf:Coupl:boundCouplProb}
 \tilde{\prob}(S) \ge 1 - (1- C_2)^{N(s,C)+1}
\end{equation}
if $\eta$ is small enough.

\medskip
Now we show how this coupling can be used to prove the proposition. First rewrite the quotient in \eqref{eq:lem:Coupling}
\begin{equation}
 \frac{\prob_{\eta}(V^{\bullet} \,|\, V)}{\prob_{\eta}(W^{\bullet} \,|\, W)}
 = \frac{\tilde{\prob}((V^{\bullet} \times \Omega^{\eta}) \cap S) + \tilde{\prob}(V^{\bullet} \times \Omega^{\eta} \,|\, S^{c})\tilde{\prob}(S^{c})}{\tilde{\prob}((\Omega^{\eta} \times W^{\bullet}) \cap S) + \tilde{\prob}(\Omega^{\eta} \times W^{\bullet} \,|\, S^{c})\tilde{\prob}(S^{c})}.\label{eq:pf:Coupl:mainQuo}
\end{equation}
We claim that
\begin{eqnarray}
 \tilde{\prob}(V^{\bullet} \times \Omega^{\eta} \,|\, S^{c}) & \asymp & \prob_{\eta}(z \lra \HP \setminus B_{2s}(z)); \label{eq:pf:Coupl:Bnd:1}\\
 \tilde{\prob}(V^{\bullet} \times \Omega^{\eta} \,|\, S) & \asymp & \prob_{\eta}(z \lra \HP \setminus B_{2s}(z)); \label{eq:pf:Coupl:Bnd:2}
\end{eqnarray}
for $\eta$ small enough. Similarly for $\Omega^{\eta} \times W^{\bullet}$.
Applying these claims together with \eqref{eq:pf:Coupl:whenCoupled} and the fact that $\tilde{\prob}(S^{c})$ converges to
zero as $s$ tends to zero, uniformly in $\eta$ as follows from \eqref{eq:pf:Coupl:boundCouplProb},
proves the proposition.

It remains to prove the claims \eqref{eq:pf:Coupl:Bnd:1} and \eqref{eq:pf:Coupl:Bnd:2}.
At first sight one might think that these bounds are easy consequences of RSW, FKG arguments.
This is not completely true since we have to deal with the condition that the coupling was not successful,
respectively successful, which are neither increasing nor decreasing events.
Recall the definition of $RM$ in \eqref{eq:pf:Coupl:defRM}. Let $PN := \HP \setminus RM$.
It is sufficient to show that, for all suitable $\nu_{PN} \times \omega_{PN}$,
\begin{equation}\label{eq:pf:Coupl:Bnd:1:claim}
 \tilde{\prob}(V^{\bullet} \times \Omega^{\eta} \,|\, \nu_{PN} \times \omega_{PN}) \asymp \prob_{\eta}(z \lra \HP \setminus B_{2s}(z)).
\end{equation}
First note that it follows from the coupling procedure that
\begin{equation*}
 \tilde{\prob}(V^{\bullet} \times \Omega^{\eta} \,|\, \nu_{PN} \times \omega_{PN}) = \prob_{\eta}(V^{\bullet} \,|\, V \cap \{\nu_{PN}\}).
\end{equation*}

\medskip
First we prove that in \eqref{eq:pf:Coupl:Bnd:1:claim}, the left hand side is less than or equal to a constant times the right hand side.
To do this we introduce the event $B$, that there is an open (semi-)circuit in $A(z;s,2s)$.
We will prove this upper bound by showing that there exist universal constants $C_3,C_4 >0$ such that, for all suitable $\nu_{PN}$
\begin{eqnarray}
 \prob_{\eta}(V^{\bullet}\cap B \,|\, V \cap \{\nu_{PN}\}) & \ge & C_3\, \prob_{\eta}(V^{\bullet} \,|\,V \cap \{ \nu_{PN}\}); \label{eq:pf:Coupl:Bnd:1:claim:lower}\\
 \prob_{\eta}(V^{\bullet}\cap B \,|\,V \cap \{\nu_{PN}\}) & \le & C_4\, \prob_{\eta}(z \lra \HP \setminus B_{2s}(z)). \label{eq:pf:Coupl:Bnd:1:claim:upper}
\end{eqnarray}
First we consider the lower bound \eqref{eq:pf:Coupl:Bnd:1:claim:lower}.
Let $\nu_{PN}$ be arbitrary.
Using Lemma \ref{lem:genFKG} and standard RSW, FKG arguments we get that
\begin{eqnarray*}
 \prob_{\eta}(V^{\bullet}\cap B \,|\, V\cap \{\nu_{PN}\}) & \ge & \prob_{\eta}(B)\prob_{\eta}(V^{\bullet}\,|\, V\cap \{\nu_{PN}\}).\nonumber\\
 & \ge & C_3\, \prob_{\eta}(V^{\bullet}\,|\, V\cap \{\nu_{PN}\})\nonumber
\end{eqnarray*}
This proves \eqref{eq:pf:Coupl:Bnd:1:claim:lower}.

Next we prove the upper bound \eqref{eq:pf:Coupl:Bnd:1:claim:upper}.
Therefore let $\Gamma$ denote the outermost open (semi-)circuit in
$A(z;s,2s)$.
Since $V$ is an $(s,C)$-one-arm like event, we have by Definition \ref{defi:OneArmLike},
\begin{equation}\label{eq:pf:Coupl:bnd:1:subset}
 \bigcup_{\gamma} V_{ext(\gamma)} \cap \{\Gamma = \gamma\} \cap \{I(z,s) \lra \gamma\} \,\,\subset\,\, V.
\end{equation}
This, together with standard RSW, FKG arguments, implies that there exists a constant $C_5>0$ such that
\begin{eqnarray}
 \prob_{\eta}(B\cap V \,|\, \nu_{PN}) & \ge & \prob_{\eta}(B\cap V_{ext(\Gamma)} \cap \{I(z,s) \lra \Gamma\} \,|\, \nu_{PN})\nonumber\\
 & \ge & C_5\, \prob_{\eta}(B\cap V_{ext(\Gamma)} \,|\, \nu_{PN}),\label{eq:pf:Coupl:VVext}
\end{eqnarray}
since $\prob_{\eta}(I(z,s) \lra \Gamma \,|\, B;\, V_{ext(\Gamma)};\, \nu_{PN}) \ge C_5$.
Hence
\begin{eqnarray}
 \prob_{\eta}(V^{\bullet}\cap B \,|\,V\cap \{\nu_{PN}\}) & \le & \prob_{\eta}(\{z \lra \Gamma\} \cap B \,|\,V\cap \{\nu_{PN}\}) \nonumber\\
 & \le & \prob_{\eta}(z \lra \HP \setminus B_{s}(z))\cdot \frac{\prob_{\eta}(B\cap V_{ext(\Gamma)} \,|\, \nu_{PN})}{\prob_{\eta}(V \,|\, \nu_{PN})} \nonumber\\
 & \le & \frac{1}{C_5 C_6} \prob_{\eta}(z \lra \HP \setminus B_{2s}(z))\cdot \frac{\prob_{\eta}(B\cap V \,|\, \nu_{PN})}{\prob_{\eta}(V \,|\, \nu_{PN})} \nonumber\\
 & \le & \frac{1}{C_5 C_6} \prob_{\eta}(z \lra \HP \setminus B_{2s}(z)),\label{eq:pf:Coupl:bnd:1:part:2}
\end{eqnarray}
where we used in the first inequality Definition \ref{defi:pointVersion}. In the second inequality
we used the fact that
$V \subset V_{ext(\Gamma)}$ together with the fact that $\{ z \lra \Gamma\}$
is independent of everything outside $\Gamma$ (which exists because of $B$).
The third inequality follows from \eqref{eq:pf:Coupl:VVext} and the existence of a universal constant $C_6>0$ such that\newline
$\prob_{\eta}(z \lra \HP \setminus B_{2s}(z)) \ge C_6\, \prob_{\eta}(z \lra \HP \setminus B_{s}(z))$.
This gives the desired inequality \eqref{eq:pf:Coupl:Bnd:1:claim:upper} and completes the proof of
the upper bound in \eqref{eq:pf:Coupl:Bnd:1:claim}.

\medskip
Next we consider the lower bound in \eqref{eq:pf:Coupl:Bnd:1:claim}. We prove that
\begin{equation}\label{eq:pf:Coupl:Bnd:2:claim}
 \prob_{\eta}(V^{\bullet} \,|\,V\cap \{\nu_{PN}\}) \ge C_3\, \prob_{\eta}(z \lra \HP \setminus B_{2s}(z)).
\end{equation}
To prove this, we again use the event $B$. The inequality \eqref{eq:pf:Coupl:Bnd:2:claim}
follows immediately from the following inequality
\begin{equation}\label{eq:pf:Coupl:Bnd:2:claim:lower}
 \prob_{\eta}(V^{\bullet}\cap B \,|\,V \cap \{\nu_{PN}\}) \ge C_3\, \prob_{\eta}(z \lra \HP \setminus B_{2s}(z)), 
\end{equation}
where $C_3>0$ is the same as in \eqref{eq:pf:Coupl:Bnd:1:claim:lower}.
Similarly to \eqref{eq:pf:Coupl:bnd:1:subset}, but now using Definition \ref{defi:pointVersion}, we have
\begin{equation}\label{eq:pf:Coupl:bnd:2:subset}
 \bigcup_{\gamma} \{\Gamma = \gamma\} \cap V_{ext(\gamma)} \cap \{z \lra \gamma\} \,\,\subset\,\, V^{\bullet},
\end{equation}
where $\Gamma$ is the outermost circuit in $A(z;s,2s)$.
Hence
\begin{eqnarray}
 \prob_{\eta}(V^{\bullet}\cap B \,|\, V\cap \{\nu_{PN}\}) & \stackrel{\eqref{eq:pf:Coupl:bnd:2:subset}}{\ge} & \sum_{\gamma} \frac{\prob_{\eta}(\{\Gamma = \gamma\} \cap V_{ext(\gamma)} \cap \{z \lra \gamma\} \cap \{\nu_{PN}\})}{\prob_{\eta}(V\cap \{\nu_{PN}\})},\nonumber\\
 & \ge & \prob_{\eta}(z \lra \HP \setminus B_{2s}(z)) \sum_{\gamma} \frac{\prob_{\eta}(\{\Gamma = \gamma\} \cap V_{ext(\gamma)} \cap \{\nu_{PN}\})}{\prob_{\eta}(V\cap \{\nu_{PN}\})},\nonumber\\
 & \ge & \prob_{\eta}(z \lra \HP \setminus B_{2s}(z)) \frac{\prob_{\eta}(B \cap V \cap \{\nu_{PN}\})}{\prob_{\eta}(V\cap \{\nu_{PN}\})}.
\end{eqnarray}
It follows from Lemma \ref{lem:genFKG} together with the fact that $\prob_{\eta}(B) \ge C_3$ that
\begin{equation}
\prob_{\eta}(B \cap V \cap \{\nu_{PN}\}) \ge C_3\cdot\prob_{\eta}(V \cap \{\nu_{PN}\}).
\end{equation}
This completes the proof of \eqref{eq:pf:Coupl:Bnd:2:claim:lower}
and finishes the proof of Proposition \ref{lem:Coupling}\qed


\subsection{Proof of Theorem \ref{thm:mainTheorem:factPPP}.}\label{subsec:pfOfThmfactPPP}
Let $u_1,u_2,w$ be fixed. Because of Theorem \ref{thm:factFormBI}
it is sufficient to show that for every $\varepsilon > 0$, there exists
$s > 0$, such that $\forall s_1,s_2,s_3 < s:$ $\exists \eta_{0}>0$ with the property that
\begin{equation}\label{eq:pf:mainThm:mainStep}
 \left| \frac{\prob_{\eta}(u_1 \lra u_2 \lra w \,|\, E_{u_1,u_2,w}^{s_1,s_2,s_3})^{2}}{\prob_{\eta}(u_1 \lra u_2 \,|\, E_{u_1,u_2}^{s_1,s_2})\,\prob_{\eta}(u_1 \lra w \,|\, E_{u_1,w}^{s_1,s_3})\,\prob_{\eta}(u_2 \lra w \,|\, E_{u_2,w}^{s_2,s_3})} - 1\right| < \varepsilon,
\end{equation}
for all $\eta < \eta_{0}$.

In order to prove \eqref{eq:pf:mainThm:mainStep} we define the following events:
\begin{eqnarray}
 E_{u_1,u_2}^{s_1,\bullet} & := & \{[u_1, u_1 +s_1] \cap \cl(u_2) \neq \emptyset\}; \label{eq:defi:E2}\\
 E_{u_1,w}^{\bullet,s_3} & := & \{ \rho(w,\cl(u_1)) < s_3 \}; \nonumber\\
 E_{u_2,w}^{\bullet,s_3} & := & \{ \rho(w,\cl(u_2)) < s_3 \}; \nonumber\\
 E_{u_1,u_2,w}^{s_1,\bullet,s_3} & := & \{ [u_1, u_1 +s_1] \cap \cl(u_2) \neq \emptyset \} \cap \{ \rho(w,\cl([u_1,u_1+s_1])) < s_3 \}; \nonumber\\
 E_{u_1,u_2,w}^{\bullet,\bullet,s_3} & := & \{ u_1 \lra u_2 \} \cap \{ \rho(w,\cl(u_1)) < s_3 \}.\nonumber
\end{eqnarray}
Let $C := (\min\{ |u_1 - u_2|, |u_1 - w|, |u_2 - w|\})/(2\sqrt{2})$.
We claim the following about the events defined in \eqref{eq:defi:E1} and \eqref{eq:defi:E2}.
\begin{enumerate}
 \item  Every event of the form $E_{a_1,a_2,a_3}^{s_1,s_2,s_3}$ or $E_{a_1,a_2}^{s_1,s_2}$
where the $a_{i}$'s are in $\{ u_1,u_2,w \}$ and each $s_{i}$ is in $\RR_{+}$ or $s_i = \bullet$,
defined in \eqref{eq:defi:E1} and \eqref{eq:defi:E2},
is, for each $s_j \neq \bullet$ an $(s_j,C)$-one-arm like event around $a_j$. For example
$E_{u_1,u_2,w}^{s_1,\bullet,s_3}$ is an $(s_1,C)$-one-arm like event around $u_1$,
and an $(s_3,C)$-one-arm like event around $w$.
 \item The events $\{u_1 \lra u_2\}$, $\{u_1 \lra w\}$, $\{u_2 \lra w\}$, $\{u_1 \lra u_2 \lra w\}$
are point versions of respectively  $E_{u_1,u_2}^{s_1,\bullet}$, $E_{u_1,w}^{\bullet,s_3}$,
$E_{u_2,w}^{\bullet,s_3}$ and $E_{u_1,u_2,w}^{\bullet,\bullet,s_3}$.
 \item Each event in \eqref{eq:defi:E2} is a point version of the corresponding event
 $E_{a_1,a_2,a_3}^{s_1,s_2,s_3}$ or $E_{a_1,a_2}^{s_1,s_2}$, where the "$\bullet$" is replaced by a positive number $s_j$.
 E.g. $E_{u_1,u_2,w}^{\bullet,\bullet,s_3}$ is a point version of $E_{u_1,u_2,w}^{s_1,\bullet,s_3}$ and
 $E_{u_1,u_2}^{s_1,\bullet}$ is a point version of $E_{u_1,u_2}^{s_1,s_2}$.
 \item Each pair of events of the form $E_{a_1,a_2,a_3}^{s_1,s_2,s_3}$ and $E_{a_1,a_2}^{s_1,s_2}$
where the $a_{i}$'s are in $\{ u_1,u_2,w \}$ and each $s_{i}$ is in $\RR_{+}$ or $s_i = \bullet$,
defined in \eqref{eq:defi:E1} and \eqref{eq:defi:E2},
are, for each $j$ where both events have $s_j \neq \bullet$, $(s_j,C)$-comparable around $a_j$.
For example the events
$E_{u_1,u_2}^{s_1,s_2}, E_{u_1,w}^{s_1,s_3}, E_{u_1,u_2,w}^{s_1,s_2,s_3}, E_{u_1,u_2}^{s_1,\bullet}, E_{u_1,u_2,w}^{s_1,\bullet,s_3}$ are
pairwise $(s_1,C)$-comparable around $u_1$.
\end{enumerate}


Before we give proofs of these claims we show how Theorem \ref{thm:mainTheorem:factPPP} follows from them.
We factorize the numerator in \eqref{eq:pf:mainThm:mainStep} as follows
\begin{eqnarray}\label{eq:pf:mainThm:numerator}
 \lefteqn{\prob_{\eta}(u_1 \lra u_2 \lra w\,|\, E_{u_1,u_2,w}^{s_1,s_2,s_3})^{2}} \\
 & = & \prob_{\eta}(u_1 \lra u_2 \lra w \,|\, E_{u_1,u_2,w}^{\bullet,\bullet,s_3})^{2}\cdot\prob_{\eta}(E_{u_1,u_2,w}^{\bullet,\bullet,s_3} \,|\, E_{u_1,u_2,w}^{s_1,\bullet,s_3})^{2}\cdot\prob_{\eta}(E_{u_1,u_2,w}^{s_1,\bullet,s_3}\nonumber \,|\, E_{u_1,u_2,w}^{s_1,s_2,s_3})^{2}.
\end{eqnarray}
The probabilities in the denominator in \eqref{eq:pf:mainThm:mainStep} can be factorized as follows
\begin{eqnarray}
 \prob_{\eta}(u_1 \lra u_2 \,|\, E_{u_1,u_2}^{s_1,s_2}) & = & \prob_{\eta}(u_1 \lra u_2 \,|\, E_{u_1,u_2}^{s_1,\bullet})\prob_{\eta}(E_{u_1,u_2}^{s_1,\bullet} \,|\, E_{u_1,u_2}^{s_1,s_2})\label{eq:pf:mainThm:den1}\\
 \prob_{\eta}(u_1 \lra w \,|\, E_{u_1,w}^{s_1,s_3}) & = & \prob_{\eta}(u_1 \lra w \,|\, E_{u_1,w}^{\bullet,s_3})\prob_{\eta}(E_{u_1,w}^{\bullet,s_3} \,|\, E_{u_1,w}^{s_1,s_3})\label{eq:pf:mainThm:den2}\\
 \prob_{\eta}(u_2 \lra w \,|\, E_{u_2,w}^{s_2,s_3}) & = & \prob_{\eta}(u_2 \lra w \,|\, E_{u_2,w}^{\bullet,s_3})\prob_{\eta}(E_{u_2,w}^{\bullet,s_3} \,|\, E_{u_2,w}^{s_2,s_3})\label{eq:pf:mainThm:den3}.
\end{eqnarray}
Plugging this into the quotient in \eqref{eq:pf:mainThm:mainStep} and applying Proposition \ref{lem:Coupling} to
the 6 pairs of $(s_i,C)$-comparable events completes the proof.

\medskip It remains to prove claims 1-4 above.
Some of these claims follow immediately, for the others we use two standard properties of conformal radius.
The first is \eqref{eq:pf:mainThm:equivConfRad}.
The second property is \emph{monotonicity}: the conformal radius is non-decreasing as the domain $A$ decreases,
(as is well known and follows easily from Schwarz' Lemma. See for example \cite{A73Book}).

We prove claim 1 for a particular event, namely $E_{u_1,w}^{\bullet,s_3}$.\newline
\textbf{(a)} It is increasing: Let $\omega \in E_{u_1,w}^{\bullet,s_3}$ and
$\nu \ge \omega$, then $\cl(u_1)(\omega) \subset \cl(u_1)(\nu)$. Here $\cl(u_1)(\omega)$ means
the cluster of $u_1$ under the configuration $\omega$.
Thus by monotonicity of the conformal radius $\rho(w,\cl(u_1)(\nu)) \le \rho(w,\cl(u_1)(\omega)) < s_3$
and $\nu \in E_{u_1,w}^{\bullet,s_3}$.\newline
\textbf{(b)} $E_{u_1,w}^{\bullet,s_3} \subset \{ B_{s_3}(w) \lra \HP \setminus B_{C}(w) \}$: Suppose that
$\omega \in E_{u_1,w}^{\bullet,s_3}$. It follows from \eqref{eq:pf:mainThm:equivConfRad} that
$\min_{x \in \cl(u_1)} |w-x| < s_3$. Further $\sqrt{2}C \le |u_1-w|/2$, which implies that
$\omega \in \{ B_{s_3}(w) \lra \HP \setminus B_{C}(w) \}$.\newline
Let $\gamma$ be an arbitrary (semi-)circuit in $A(w;s_3,C)$. Let $D := E_{u_1,w}^{\bullet,s_3}$\newline
\textbf{(c)} $\{ \gamma \textrm{ open}\} \cap D_{ext(\gamma)} \cap D_{int(\gamma)} \subset D$: Let $\omega \in D_{int(\gamma)}$
and $\nu \in D_{ext(\gamma)}$. By definition there exists $\tilde{\nu}$ such that $\nu_{ext(\gamma)} \times \tilde{\nu} \in D$.
With the second inequality in \eqref{eq:pf:mainThm:equivConfRad} this implies that $u_1 \lra \gamma$ in $ext(\gamma)$.
Next let $\tilde{\omega}$ be such that $\omega_{int(\gamma)} \times \tilde{\omega} \in D$.
Then it is easy to see that $\cl(u_1)(\omega_{int(\gamma)} \times \tilde{\omega}) \cap int(\gamma) \subset \cl(\gamma)(\omega) \cap int(\gamma)$.
Monotonicity of the conformal radius implies now that
\[
 \rho\left(w,\cl(\gamma)(\omega)\right) \le \rho\left(w,\cl(u_1)(\omega_{int(\gamma)} \times \tilde{\omega})\right) < s_3
\]
Let $\upsilon := \omega_{int(\gamma)} \times \{1\}^{\gamma} \times \nu_{ext(\gamma)}$.
Note that $\cl(u_1)(\upsilon) \cap int(\gamma) = \cl(\gamma)(\omega) \cap int(\gamma)$. Thus
$\rho(w,\cl(u_1)(\upsilon)) = \rho(w,\cl(\gamma)(\omega))$, and hence $\upsilon \in D$.\newline
\textbf{(d)} $\{I(w,s_3) \lra \gamma\} \subset D_{int(\gamma)}$: Let $\omega \in \{I(w,s_3) \lra \gamma\}$
and $\nu \in \{ u_1 \lra \gamma \}$. Then the first inequality in \eqref{eq:pf:mainThm:equivConfRad} implies that
$\omega_{int(\gamma)} \times \{1\}^{\gamma} \times \nu_{ext(\gamma)} \in D$, hence $\omega \in D_{int(\gamma)}$.\newline
This completes the proof of claim 1 for this particular event. The proofs for the other events and claims
are very similar and we omit them. \qed

\subsection{Proof of Theorem \ref{thm:mainTheorem:factPPI}.}\label{subsec:pfOfThmfactPPI}
We will use the notation
\begin{equation}
 E_{u_1,u_2,w}^{s_1,\bullet,\bullet} := \left\{\{u_2, w\} \subset \cl([u_1,u_1+s])\right\}.
\end{equation}
With this notation we can write the quotient in \eqref{eq:thm:PPI} as
\begin{equation}
\frac{\prob(\{u_2, w\} \subset \cl([u_1,u_1+s]))}{\prob(w \in \cl([u_1,u_1+s]) )\,\, \prob(u_2 \in \cl([u_1,u_1+s]))} = \frac{\prob_{\eta}(E_{u_1,u_2,w}^{s,\bullet,\bullet})}{\prob_{\eta}(E_{u_1,w}^{s,\bullet})\prob_{\eta}(E_{u_1,u_2}^{s,\bullet})}.\label{eq:pf:thmPPI:shortStat}
\end{equation}
Similarly to the proof of Theorem \ref{thm:mainTheorem:factPPP} we factorize this as
follows
\begin{eqnarray}
\lefteqn{\frac{\prob_{\eta}(E_{u_1,u_2,w}^{s,\bullet,\bullet})}{\prob_{\eta}(E_{u_1,w}^{s,\bullet})\prob_{\eta}(E_{u_1,u_2}^{s,\bullet})}}\label{eq:pf:thmPPI:fact}\\
 & = & \frac{\prob_{\eta}(E_{u_1,u_2,w}^{s,\bullet,\bullet} \,|\, E_{u_1,u_2,w}^{s,\bullet,s_3})}{\prob_{\eta}(E_{u_1,w}^{s,\bullet} \,|\, E_{u_1,w}^{s,s_3})}\cdot\frac{\prob_{\eta}(E_{u_1,u_2,w}^{s,\bullet,s_3} \,|\, E_{u_1,u_2,w}^{s,s_2,s_3})}{\prob_{\eta}(E_{u_1,u_2}^{s,\bullet} \,|\, E_{u_1,u_2}^{s,s_2})}\cdot\frac{\prob_{\eta}(E_{u_1,u_2,w}^{s,s_2,s_3})}{\prob_{\eta}(E_{u_1,w}^{s,s_3})\prob_{\eta}(E_{u_1,u_2}^{s,s_2})}.\nonumber
\end{eqnarray}
The first two ratio's converge to 1 by Proposition \ref{lem:Coupling}, uniformly in $\eta$. Namely
the involved events are point versions and $(s,C)$-comparable, by similar arguments as in
the proof of Theorem \ref{thm:mainTheorem:factPPP}. We claim that the ratio
\begin{equation}
 \frac{\prob_{\eta}(E_{u_1,u_2,w}^{s,s_2,s_3})}{\prob_{\eta}(E_{u_1,w}^{s,s_3})\prob_{\eta}(E_{u_1,u_2}^{s,s_2})}
\end{equation}
converges to the function $\psi(u_1,s,u_2,w)$, as $\eta, s_2,s_3$ tend to zero.
To prove this claim we note that
\begin{equation}\label{eq:pf:thmPPI:claim}
 \frac{\prob_{\eta}(E_{u_1,u_2,w}^{s,s_2,s_3})}{\prob_{\eta}(E_{u_1,w}^{s,s_3})\prob_{\eta}(E_{u_1,u_2}^{s,s_2})}
 = \frac{s_3^{-5/48}\cdot \prob_{\eta}(E_{u_1,u_2,w}^{s,s_2,s_3} \,|\, E_{u_1,u_2}^{s,s_2})}{s_3^{-5/48}\cdot\prob_{\eta}(E_{u_1,w}^{s,s_3})}.
\end{equation}
Theorem \ref{thm:condProbBI} and Lemma \ref{lem:poiInt} imply that the following limit of \eqref{eq:pf:thmPPI:claim} exists:
First send $\eta$ to zero, after that send $s_2$ to zero and finally let $s_3$ go to zero.
This, together with the uniform convergence in $\eta$ of the first two ratio's in 
\eqref{eq:pf:thmPPI:fact}, implies that the limit in \eqref{eq:thm:PPI} exists and is equal to
\begin{equation}\label{eq:pf:ThmPPI:Psi:part:1}
 \frac{\pi^{5/48}|\Psi'_{u_1,s,u_2}(w)|^{5/48}G\left(\Re (\Psi_{u_1,s,u_2}(w)), \Im (\Psi_{u_1,s,u_2}(w))\right)}{2^{5/48}H(0)\cdot |\phi'(w)|^{5/48}(\sin(\pi\omega/2))^{1/3}},
\end{equation}
where $\Psi, G, \phi, H,\omega$ are as in Theorem \ref{thm:condProbBI} and Lemma \ref{lem:poiInt}.

To finish the proof of Theorem \ref{thm:mainTheorem:factPPI} we have to simplify \eqref{eq:pf:ThmPPI:Psi:part:1}
and show that it is equal to the function $\psi(u_1,s,u_2,w)$ given in that Theorem.
Hereto let $\Pi:\HP \to \HP$ be a conformal map such that the points $u_1,u_1+s,u_2$ are mapped to $-1,1,\infty$ respectively.
Let $\tilde{w} = \Pi(w)$.
Let $\tilde{\Psi}:\HP \to \mathbb{S}$ be the conformal map, such that $\Psi = \tilde{\Psi} \circ \Pi$,
thus
\[
 \tilde{\Psi}(z) = \frac{-\imag}{\pi}\arcsin(z)+\frac{1}{2}\imag.
\]
Further let $\tilde{\phi}$ be the conformal map such that $\phi = \tilde{\phi} \circ \Pi$.
We have that
\begin{equation}\label{eq:pf:ThmPPI:derivatives}
 |\phi'(w)| = \frac{|\Pi'(w)|}{2\Im(\tilde{w})}, \qquad 
 |\Psi'(w)| = \frac{|\Pi'(w)|}{\pi\sqrt{|1-\tilde{w}^2|}}.
\end{equation}
Recall that $x=\Re (\Psi_{u_1,s,u_2}(w)),\, y= \Im (\Psi_{u_1,s,u_2}(w))$
and $\Psi_{u_1,s,u_2}(w) = \tilde{\Psi}(\tilde{w})$, thus 
\begin{eqnarray*}
 \sinh(\pi x) & = & \sinh(\Im(\arcsin(\tilde{w}))),\\
 \sin(\pi y) & = & \cos(\Re(\arcsin(\tilde{w}))).
\end{eqnarray*}
It follows from standard formulas for hyperbolic functions that
\begin{eqnarray}
 \sinh(\pi x)^{2}\sin(\pi y)^{2} & = & \Im(\tilde{w})^2, \label{eq:pf:ThmPPI:part:2}\\
 \sinh(\pi x)^{2}+\sin(\pi y)^{2} & = & |1-\tilde{w}^2|.\label{eq:pf:ThmPPI:part:3}
\end{eqnarray}
Further note that
\begin{eqnarray}
 \lefteqn{\left(\frac{1}{\sinh(\pi x)}\right)^{1/3}\left(\frac{\sinh(\pi x)^{2}\sin(\pi y)^{2}}{\sinh(\pi x)^{2}+\sin(\pi y)^{2}} \right)^{11/96}} \label{eq:pf:ThmPPI:part:4}\\
 & = & \left(\frac{\sin(\pi y)^{2}}{\sinh(\pi x)^{2}+\sin(\pi y)^{2}} \right)^{1/6}\left(\frac{\sinh(\pi x)^{2}+\sin(\pi y)^{2}}{\sinh(\pi x)^{2}\sin(\pi y)^{2}} \right)^{5/96}.\nonumber
\end{eqnarray}
Putting together the definition of $G$ in \eqref{eq:thmCondProbBI:DefG} and equations \eqref{eq:pf:ThmPPI:derivatives} - \eqref{eq:pf:ThmPPI:part:4} gives that \eqref{eq:pf:ThmPPI:Psi:part:1} is equal to
\begin{equation}\label{eq:pf:ThmPPI:part:5}
 \frac{e^{\pi x/3}H(x)}{H(0)}\cdot \left(\frac{\cos(\Re(\arcsin(\tilde{w})))}{\sqrt{|1-\tilde{w}^2|}\sin(\pi\omega/2)}\right)^{1/3}.
\end{equation}
Recall that $\omega\pi$ is equal to the angle at $\tilde{w}$ in the triangle with corners $-1,1,\tilde{w}$.
It follows easily that
\begin{equation*}
 \sin(\pi \omega/2) = \sqrt{\frac{1}{2} - \frac{|\tilde{w}|^2-1}{2|1-\tilde{w}^2|}},
\end{equation*}
and from formulas for hyperbolic functions, including \eqref{eq:pf:ThmPPI:part:3}, that
\begin{equation*}
 2\cos(\Re(\arcsin(\tilde{w})))^2 = |1-\tilde{w}^2|+1-|\tilde{w}|^2,
\end{equation*}
which together imply that the last factor in \eqref{eq:pf:ThmPPI:part:5} equals 1.
This completes the proof of Theorem \ref{thm:mainTheorem:factPPI}.\qed


\medskip\noindent\textbf{Acknowledgments.} The author would like to thank Rob van den Berg for stimulating discussions and
comments on earlier drafts of this paper.

\bibliographystyle{elsarticle-num}
\bibliography{mybiblio}

\end{document}